\documentclass[11pt]{article}
\usepackage{amsfonts}
\usepackage{amsthm}
\usepackage{color}
\usepackage{graphicx}

\newtheorem{theorem}{Theorem}[section]
\newtheorem{lemma}{Lemma}[section]

\newtheorem{corollary}{Corollary}[section]

\theoremstyle{definition}
\newtheorem{definition}{Definition}[section]

\newcommand\noi{\noindent}
\newcommand\ba {{\bf{a}}}
\newcommand\bb {{\bf{b}}}

\newcommand\be {{\bf{e}}}

\newcommand\bh {{\bf{h}}}

\newcommand\bt {{\bf{t}}}

\newcommand\bv {{\bf{v}}}
\newcommand\bw {{\bf{w}}}
\newcommand\bx {{\bf{x}}}
\newcommand\by {{\bf{y}}}

\newcommand\bA {{\bf{A}}}

\newcommand\bC {{\bf{C}}}

\newcommand\bU {{\bf{U}}}
\newcommand\bV {{\bf{V}}}
\newcommand\bW {{\bf{W}}}
\newcommand\bX {{\bf{X}}}
\newcommand\bY {{\bf{Y}}}
\newcommand\bZ {{\bf{Z}}}

\newcommand\identidad {\mbox{\bf I}}

\newcommand\bmu {\mbox{\boldmath $\mu$}}

\newcommand\bLam {\mbox{\boldmath $\Lambda$}}
\newcommand\bGa {\mbox{\boldmath $\Gamma$}}

\newcommand\bSi {\mbox{\boldmath $\Sigma$}}

\newcommand{\trasp}{^{\mbox{\sc t}}}

\def\argmin{\mathop{\rm argmin}}

\newcommand{\igualdist}{ \buildrel{{\cal D}}\over \sim}

\def\argmin{\mathop{\rm argmin}}
\def\square{\ifmmode\sqr\else{$\sqr$}\fi}

\def\dst{\displaystyle}
\def\noi{\noindent}
\def\new{\newline}

\def\real{\hbox{$\displaystyle I\hskip -3pt R$}}
\def\realito{\hbox{\footnotesize $\mathbb{R}$}}

\def\natu{\hbox{$\displaystyle I\hskip -3pt N$}}
\def\natuito{\hbox{\footnotesize$\displaystyle I\hskip -3pt N$}}

\def\uno{\mathbb{I}}
\def\square{\ifmmode\sqr\else{$\sqr$}\fi}
\def\sqr{\vcenter{
         \hrule height.1mm
         \hbox{\vrule width.1mm height2.2mm\kern2.18mm
\vrule width.1mm}
         \hrule height.1mm}}

\newcommand{\elemAleat}[1]{{#1}}

\begin{document}

\title{{Principal points and elliptical distributions
from the multivariate setting to the functional case}\thanks{ This research was partially supported by Grants X-018 from the Universidad de Buenos Aires, \sc pid \rm  5505 from \textsc{conicet} and \textsc{pav} 120 and \sc pict \rm  21407 from \textsc{anpcyt}, Argentina. The research of Lucas Bali was supported by a scholarship  of the Agencia Nacional de Promoci\'on Cient\'\i fica y Tecnol\'ogica.}}
\author{
{Juan Lucas Bali}\\
Universidad de Buenos Aires and FONCYT, Argentina\\
{and}\\
{ Graciela Boente }\\
        Universidad de Buenos Aires and CONICET, Argentina \\
}
\date{}
\maketitle

\begin{abstract}
The $k$ principal points of a random vector $\bX$ are defined as a set of points which minimize the
expected squared distance between $\bX$ and the nearest point in
the set. They are thoroughly studied in Flury (\cite{flury90}, \cite{flury93}), Tarpey \cite{tarpey95}   and Tarpey, Li and Flury \cite{tarpeyliflury}. For their treatment, the examination
is usually restricted to the family of elliptical distributions. In this paper, we present an extension of the previous results to the functional elliptical distribution case, i.e., when dealing with random elements over a separable Hilbert space ${\cal H}$. Principal points for gaussian processes were defined in Tarpey and Kinateder \cite{tarpeykini}. In this paper, we generalize the concepts of principal points, self-consistent points and elliptical distributions so as to fit them in this functional framework. Results linking self-consistency and the eigenvectors of the covariance operator are re-obtained in this new setting as well as an explicit formula for the $k=2$ case so as to include elliptically distributed random elements in ${\cal H}$.
\end{abstract}

\newline{\em Key words and phrases:}

\vskip0.3in
\noi \textbf{Corresponding Author}
\new Graciela Boente
\new Moldes 1855, 3$^o$ A
\new Buenos Aires, C1428CRA, Argentina
\new email: gboente@dm.uba.ar
\normalsize

\newpage
\section{Introduction}

\subsection{Motivation}

Inside statistics there exists lots of situations where the collected data may not be represented with classic schemes like numbers or numeric vectors and so, sometimes a functional representation is more appropriate. For example, consider results of an electrocardiogram (EGC) or the study of the temperature in a weather station, which lend themselves to this new framework (see, for instance,  Ramsay  and Silverman \cite{ramsay97functional} for more examples). A classical discretization of the data as a sequence of numbers may lose some functional characteristics like smoothness and continuity. For this reason, in the last decades different methods appeared to handle this new kind of data. In an informal way, we may say that a functional datum is a random variable (element would be a better word) that takes its values in a functional space, instead of a finite dimensional one. In this paper, we will study some
fundamental concepts of this development, which in a way will result in a mixture between statistics and functional analysis over  Hilbert spaces.
The main idea is to mix together, in a very general family of distributions, some notions of principal components and principal points. In a multivariate setting, those developments were mainly done by Flury and Tarpey (\cite{flury90} to \cite{ft99}, \cite{tarpey94} to \cite{tarpey992}, \cite{tarpeyflury96} and \cite{tarpeyliflury}) at the beginning of the 90's. The idea here is to adapt the results obtained therein to the functional elliptical distribution case.

Maybe the final conclusion of this work is not only the theoretical result obtained. Rather, as it was done  previously, our results show  about the possibility of doing, with some technical difficulties but not critical ones, an interesting generalization of the classical results from multivariate analysis to a more general framework, so as to gain a better comprehension of the phenomenon, as it often tends to happen when abstraction or generalization of a mathematical concept is made.

In section \ref{ellip}, we will define the notion of elliptical families. We will first remind their definition in  the multivariate case and later we will extend this definition to the functional case.  The definition of self--consistent points and principal points as well as some of their properties are stated in  section \ref{self} where we extend the results given in Flury (\cite{flury90}, \cite{flury93}), Tarpey \cite{tarpey95} and Tarpey and Flury \cite{tarpeyflury96} to the case of random elements lying in a separable Hilbert space. We also provide  a characterization that, under an hypothesis of ellipticity,  allows us to make an important link between principal components and self--consistent points. We conclude with some results that allow to compute principal points in a somewhat specific case.

\section{Elliptical families}\label{ellip}

\subsection{Review on some finite--dimensional results}\label{ellipnot}
For the sake of completeness and to fix our notation we will remind some results regarding elliptical families before extending them to the functional setting. They can be found in Muirhead \cite{muirhead:1982}, Seber \cite{seber:1984} and also in Frahm \cite{frahm}.

Let $\bX \in \mathbb{R}^d$ be a random vector. We will say that $\bX$ has an elliptical distribution, and we will denote it as $\bX \sim
{\cal E}_d(\bmu, \bSi, \phi)$, if there exists a vector $\bmu \in \mathbb{R}^d$, a positive semidefinite matrix $\bSi \in\mathbb{R}^{d\times d}$ and a function $\phi: \mathbb{R}_+ \to \mathbb{R}$ such that the characteristic function of $\bX-\bmu$ is given by $\varphi_{\bX - \bmu}(\bt) = \phi(\bt\trasp\bSi \bt)$, for all $\bt
\in \mathbb{R}^d.$ In some situations, for the sake of simplicity, we will omit the symbol $\phi$ and will denote $\bX \sim {\cal E}_d(\bmu, \bSi)$.

As it is well known, if  $\bX \sim {\cal E}_d(\bmu, \bSi, \phi)$ and $E(\bX)$  exists, then  $E(\bX)=\bmu$. Moreover, if the second order moments exist $\bSi$ is up to a constant the covariance matrix of $\bX$, i.e., $Var(\bX) = \alpha \bSi$. Even more, it is easy to see that the constant $\alpha $ equals
to $-2 \phi'(0)$, where $\phi'$ stands for the derivative of $\phi$.

The following theorem is a well known result and will also be extended in the sequel to adapt for functional random elements.

\begin{theorem}\label{elip:finito}
Let $\bX \sim {\cal E}(\bmu, \bSi, \phi)$ with $\bmu \in \mathbb{R}^d$ and $\bSi \in \mathbb{R}^{d\times d}$ a semidefinite positive matrix with $\mbox{rank}(\bSi) = r$. Let $\bmu=(\bmu_1\trasp, \bmu_2\trasp)\trasp$ and $\bX=(\bX_1\trasp, \bX_2\trasp)\trasp$ with $\bX_1$ the vector of the first $k$ coordinates of $\bX$ ($k<r$) such that $\bSi_{11}$ is not singular. Denote by
$$\bSi = \bLam \bLam\trasp = \left( \begin{array}{cc}
  \bSi_{11} & \bSi_{12} \\
  \bSi_{21} & \bSi_{22}
\end{array} \right) \in \mathbb{R}^{d\times d}
$$
with submatrixes $\bSi_{11} \in \mathbb{R}^{k\times k}$,
$\bSi_{21} \in \mathbb{R}^{(d-k)\times k}$, $\bSi_{12}=\bSi_{21}\trasp
\in \mathbb{R}^{k\times (d-k)}$ and $\bSi_{22}
\in\mathbb{R}^{(d-k)\times (r-k)}$.
Then,
\begin{enumerate}
\item[a)] $\bX\igualdist \bmu + {\cal R}\bLam \bU^{(r)}$, where $\bY\igualdist\bZ$ means that the two random vectors $\bY$ and $\bZ$ have the same distribution, $\bU^{(r)}$ is uniformly distributed over ${\cal S}^{r-1}=\{\by\in \mathbb{R}^r\;:\; \|\by\|=1\}$  and $\cal R$ and $\bU^{(r)}$ are independent.

\item[b)] Assume that the conditional random vector $\bX_2 | \bX_1 = \bx_1$ exists, then $\bX_2 | \bX_1
= \bx_1$ has an elliptical distribution ${\cal E}_{d-k}(\bmu^*, \bSi^*, \phi^*)$ where
\begin{eqnarray*}
\bmu^* &=& \bmu_2 + \bSi_{21}\bSi_{11}^{-1}(\bx_1 - \bmu_1)\\
\bSi^* &=& \bSi_{22} - \bSi_{21}\bSi_{11}^{-1}\bSi_{12},
\end{eqnarray*}
and $\phi^*$ corresponds to the characteristic generator of ${\cal R}^*\bU^{(r-k)}$ with
$${\cal R}^* \igualdist{\cal R}\sqrt{1-\beta}\;\Big|({\cal R}\sqrt{\beta}\;\bU^{(k)}=\bC_{11}^{-1}(\bx_1- \bmu_1)).$$
Here  $\bC_{11}$ stands for the Cholesky square root of $\bSi_{11}$, $\bU^{(k)}$ is uniformly distributed in ${\cal
S}^{k-1}$,  $\beta\sim \mbox{Beta}(\frac{k}{2}, \frac{r-k}{2})$ and ${\cal R}$, $\beta$, $\bU^{(k)}$ and $\bU^{(r-k)}$ are mutually
independent.
\end{enumerate}
\end{theorem}

\subsection{Functional case}\label{ellipfun}

In this section, we will extend the definition of elliptical distributions to the case of random elements on a separable Hilbert space. The definition will be based on the one given for the multivariate case.

\begin{definition}
Let $V$ be a random element in a separable Hilbert space $\cal H$. We will say that $V$ has an elliptical distribution of parameters $\mu\in {\cal H}$ and $\bGa:{\cal H} \rightarrow {\cal H}$, with $\bGa$ a self--adjoint, positive semidefinite and compact operator, and we
will denote $\bV\sim {\cal E}(\mu, \bGa)$, if for any lineal and bounded operator $A:{\cal H} \rightarrow \real^d$ (that is, such that
$\dst\sup_{\|x\|=1}\|Ax\|<\infty$) we have that $AV$ has a multivariate elliptical distribution of parameters $A\mu$ and $\bSi=A\bGa A^* $, i.e., $AV\sim {\cal E}_{d}(A\mu, \bSi)$ where $A^*: \mathbb{R}^p\rightarrow {\cal H}$ stands for the adjoint operator of $A$.
\end{definition}

The following result shows that elliptical families in Hilbert spaces are closed through linear and bounded transformations.

\begin{lemma} \label{elip:comp}
Let $V  \sim {\cal E}(\mu, \bGa)$ an elliptical random element
in ${\cal H}_1$ of parameters $\mu$ and $\bGa$ and $A:{\cal H}_1
\rightarrow {\cal H}_2$ linear and bounded. Then $AV$ is an
elliptical random element in ${\cal H}_2$ of parameters  $A\mu$ and
$A\bGa A^*$.
\end{lemma}

Lemma \ref{elip:esp} shows that both parameters, $\mu$ and $\bGa$, that characterizes the element $V$ are respectively the expectation and the covariance operator, provide they exist. Its proof can be found in the Appendix.

\begin{lemma} \label{elip:esp}
Let $V$ be a random element in a separable Hilbert space $\cal H$ such that $\bV\sim {\cal E}(\mu, \bGa)$.
\begin{enumerate}
\item[a)] If  $E(V)$ exists, then, $E(V) = \mu$.
\item[b)] If the covariance operator, $\bGa_V$,  exists then,  $\bGa_V = {\alpha}\;\bGa$, for some $\alpha\in \mathbb{R}$.
\end{enumerate}
\end{lemma}

Based on the finite dimensional results,  one way of obtaining random elliptical elements is through the following
transformation. Let $V_1$  be a gaussian element in ${\cal H}$ with zero mean and covariance operator $\bGa_{V_1}$, and let $Z$
be a random variable with distribution $G$ independent of $V_1$. Given $\mu\in {\cal H}$, define $V=\mu+Z\, V_1$. Then, $V$ has an
elliptical distribution and if $E(Z^2)$ exists $\bGa_{V}=E(Z^2)\,\bGa_{V_1}$.

\vskip0.2in

We are interested in obtaining some properties concerning the conditional distribution of elliptical families similar to those existing in the multivariate setting. Let $V$ be a random element belonging to an elliptical family of parameter $\mu$ and $\bGa$ and let us consider in ${\cal H}$ the orthonormal basis, $\{ \phi_n \}_{n \in {\cal I}} $ ($\cal I$ countable or finite)  constructed using the eigenfunctions of the operator $\bGa$ related to the eigenvalues $\lambda_1\ge \lambda_2\ge \ldots$. Given $d\in {\cal I}$ fixed, define the closed subspaces (and so Hilbert spaces)
$${\cal H}_1 = <\phi_1, \ldots, \phi_d>\qquad {\cal H}_2 = <\phi_1, \ldots, \phi_d>^{\bot}
$$
Define over these spaces the truncating projections, that is, $P_d=P_{{\cal H}_1}: {\cal H}\rightarrow {\cal H}_1$ and $P_{{\cal H}_2}: {\cal H} \rightarrow {\cal H}_2$ such that
\begin{equation}
P_{{\cal H}_1} (\phi_i) = \left\{ \begin{array}{cc}
\phi_i  & i =1,2, \ldots , d  \\
0 & i > d
\end{array} \right.
\qquad P_{{\cal H}_2} (\phi_i) = \left\{ \begin{array}{cc}
\phi_i & i > d  \\
0 & i =1,2, \ldots , d
\end{array} \right.\;.
\label{trunca}
\end{equation}
We will make a composition of $P_d$ with the natural operator that identifies ${\cal H}_1$ with $\mathbb{R}^d$. That is, we will
consider the operator $T_d: {\cal H} \rightarrow \mathbb{R}^d$ defined as
\begin{equation}
T_d (\phi_i) = \left\{ \begin{array}{cc}
\be_i & i =1,2, \ldots , d  \\
0 & i > n
\end{array} \right.\;,
\label{deftn}
\end{equation}
with $\be_1,\ldots,\be_d$ the  vectors of the canonical base of $\mathbb{R}^d$. Then, for any $x\in{\cal H}$  we have that $T_d(x)=\sum_{j=1}^d <x,\phi_j>\be_j$. We will use $T_d$ instead of $P_d$ as a projector in many situations, because its image is  $\mathbb{R}^d$ and we will call each of them truncating projectors.

Based on these projections we can construct $\elemAleat{V_1} = P_{{\cal H}_1} \elemAleat{V} \in{\cal H}_1$, $W_1=T_d V\in \mathbb{R}^d$ and $\elemAleat{V_2} =P_{{\cal H}_2} \elemAleat{V} \in {\cal H}_2$ random elements, both of them elliptical by Lemma \ref{elip:comp}.

We have essentially split the random element in two parts, one of them being finite dimensional which will allow us to define a
conditional distribution $\elemAleat{V}_2|\elemAleat{V_1}$ following the guidelines previously established.

\begin{theorem} \label{eliptico_partido}
Let ${\cal H}$ be a separable Hilbert space. Let $V$ be a random element in ${\cal H}$ with distribution
${\cal E}(\mu, \bGa)$ with finite second moments. Without loss of generality, we can assume that  $\bGa$ is the covariance operator. Assume that $\bGa$ is Hilbert--Schmidt so that $\sum_{i=1}^\infty\lambda_i<\infty$. Let $d\in {\cal I}$ fixed and consider $V_1=P_{{\cal H}_1} V$, $W_1=T_d V$ and $V_2=P_{{\cal
H}_2} V$  with $P_{{\cal H}_1}$ and $P_{{\cal H}_2}$ defined in (\ref{trunca}) and $T_d$ defined in (\ref{deftn}). Let $\lambda_1\ge \lambda_2\ge \ldots$ be the eigenvalues of $\bGa$ and assume that $\lambda_d>0$. Then,
\begin{itemize}
\item[a)] the covariance matrix of $W_1$ given by $\bSi_{W_1}=T_d\bGa T_d^*=\mbox{diag}\left(\lambda_1,\ldots,\lambda_d\right)$ is non--singular
\item[b)] $E(V_2|W_1) = \mu_2 + \bGa_{V_2, W_1}\bSi_{W_1}^{-1} (W_1 - \bmu_1) $,
\end{itemize}
where $\bGa_{V_2, W_1}$ is the covariance operator between $V_2$ and $W_1$, $\bmu_1=E(W_1)$ and $\mu_2=E(V_2)$.
\end{theorem}

\section{Self--consistent points and principal points}\label{self}
As mentioned in the Introduction  self--consistent  and principal points were studied by Flury (\cite{flury90}, \cite{flury93}), Tarpey \cite{tarpey95} and Tarpey and Flury \cite{tarpeyflury96} in the multivariate setting. Later on, Tarpey and Kinateder \cite{tarpeykini} extended their definition and properties for gaussian processes while Tarpey \textsl{et al.} \cite{tarpeypetrova} applied principal points to estimate a set
of representative longitudinal response curves from a clinical trial. The aim of this section is to extend some of the properties previously obtained to include elliptical families.

For the sake of completeness, we remind the definition of self--consistency and principal points.

\vskip0.2in
\begin{definition}
Let ${\cal W}=\{ y_1, \ldots, y_k \}$ with  $y_i \in\cal H$, $1\le i \le k$ we define the minimum distance of $V$ to the
set $\cal W$ as $d(V, \{y_1, \ldots, y_k\}) = \min_{1\leq j \leq k} \| V - y_j\|$.
\end{definition}

\vskip0.2in
The set $\cal W$  induce a partition of the space $\cal H$ determined by the domains of attraction.

\begin{definition}
Given ${\cal W}=\{ y_1, \ldots, y_k \}$, the domain of attraction ${\cal D}_j$ of $y_j$ consists in all the elements of
$\cal H$ that have $y_j$ as the closest point of $\cal W$, that is, ${\cal D}_j=\{x \in {\cal H}: \|x-y_j\|<
\|x-y_\ell\|\;, \quad \ell\ne j\}$.
\newline {For points $x\in {\cal H}$ with equal distance to two or several $y_j$, we assign them arbitrarily  to the set with lower index $j$.}
\end{definition}

\vskip0.2in
\begin{definition}
Let $V$ be a random element in ${\cal H}$ with expectation $E(V)$. A set ${\cal W}=\{ y_1, \ldots, y_k \}$ is said to be
self-consistent for $V$ if $E(V|V\in {\cal D}_j) = y_j$.

\noindent A random element $W$ is called self-consistent for $V$  if $E(V|W) = W$.
\end{definition}

\vskip0.2in
\begin{definition}
Let $V$ be a random element in ${\cal H}$ with finite second moment. The elements $\xi_1, \ldots, \xi_k$ are called principal points of $V$ if
$$D_V(k) = E(d^2(V,\{\xi_1, \ldots, \xi_k\})) = \min_{y_j \in {\cal H}} E(d^2(V,\{y_1, \ldots, y_k\}))$$
\end{definition}

\vskip0.2in
Lemma 3 in Tarpey and Kinateder  \cite{tarpeykini} establishes for $L^2({\cal I})$ functions, with ${\cal I}$ a real bounded interval,  the well--known result in multivariate analysis that the mean of a distribution lies in the convex hull of any set of self--consistent points. Moreover, Flury \cite{flury93} established that principal points of a random vector in $\mathbb{R}^p$ are self--consistent points. This result was generalized to random functions in $L^2({\cal I})$ by Tarpey and Kinateder  \cite{tarpeykini}. The same arguments allow to establish these results for any separable Hilber space $\cal H$,  we state them without proof.

 \begin{lemma}\label{lemaautoesp}
{Let $V$ be a random element of a separable Hilbert space $\cal H$ such that $E(V)$ exists}. Then,
\begin{enumerate}
\item[a)] if $\{ y_1, \ldots, y_k \}$ is a self-consistent set, then
$E(V)$ is a convex combination of $y_1, \ldots, y_k$.
\item[b)] Moreover, if $V$ has finite second moments and the set ${\cal W}=\{ \xi_1, \ldots, \xi_k \}$ is a set of principal points for $V$, then it is self-consistent.
\end{enumerate}
\end{lemma}

As a consequence of Lemma \ref{lemaautoesp}, if $k=1$ and $V$ is a random element with self-consistent set $\{ y_1 \}$, then $y_1=E(V)$. Moreover, we will have self-consistent points whenever we have principal points.


\vskip0.2in
The following result  will allows us to assume, in the sequel, that the random element $V$ has expectation 0. It also generalizes Lemma 2.2 in Tarpey \textsl{et al.} \cite{tarpeyliflury}, to the infinite--dimensional setting. Its proof is given in the Appendix.

\begin{lemma}\label{selfunitario}
Let $V$ be a random element of a separable Hilbert space $\cal H$ and define $V_2 = \nu + \rho \;UV$ with $\nu \in {\cal H}$, $\rho$ a scalar and $U:{\cal H}\to {\cal H}$ a unitary operator, i.e., surjective and isometric.  Then, we have that
\begin{itemize}
\item[a)] If  ${\cal W}= \{ y_1, \ldots, y_k \}$ is a set of $k$ self-consistent points of $V$, then ${\cal W}_2 = \{ \nu + \rho\; Uy_1, \ldots, \nu + \rho\; Uy_k \}$ is a set of $k$ self-consistent points of $V_2$.
\item[b)] If ${\cal W}= \{ y_1, \ldots, y_k \}$ is a set of $k$ principal points of $V$, then ${\cal W}_2 = \{ \nu + \rho\; Uy_1, \ldots, \nu + \rho \;Uy_k \}$ is a set of $k$ principal points of $V_2$ and $E(d^2(V_2, {\cal W}_2)) = \rho^2E(d^2(V,{\cal W}))$.
\end{itemize}
\end{lemma}

\vskip0.2in
Lemma \ref{finitedim} is analogous to Lemma 2.3 in Tarpey \textsl{et al.} \cite{tarpeyliflury}.
\begin{lemma}\label{finitedim}
Let $V$ be a random element with expectation $0$. Let $\{y_1,\ldots, y_k\}$ be a set of $k$ self--consistent points of $V$ spanning a subspace $\cal M$ of
dimension $q$, with an orthonormal basis $\{e_1, \ldots, e_q\}$. Then, the random vector of $\mathbb{R}^q$ defined by
$\bX=(X_1,\ldots,X_q)\trasp$ with $X_i=<e_i, V>$ will have ${\cal W}=\{\bw_{j}\}_{1\le j\le k}$ with $\bw_j=(w_{1j},\ldots,w_{qj})\trasp$
and $w_{ij}=<e_i, y_j>$ as self--consistent set.
\end{lemma}

\vskip0.2in
The notion of best $k-$point approximation  has been considered  by Tarpey \textsl{et al.} \cite{tarpeyliflury} for finite--dimensional random elements. It extends immediately to elements on a Hilbert space.

\begin{definition} Let $W\in {\cal H}$ be a discrete random element, jointly distributed with the random element $V\in {\cal H}$ and denote ${\cal S}(W)$ the support of $W$. The random element $W$ is a best $k-$point approximation  to $V$ if ${\cal S}(W)$ contains exactly $k$ different elements
$y_1, \ldots, y_k$ and $E(\|V - W\|^2) \leq
E(\|V-Z\|^2)$ for any random element $Z\in {\cal H}$ whose support has at most $k$ points, i.e.,
$\# {\cal S}(Z)\le k$.
\end{definition}

\vskip0.2in
The following result is the infinite--dimensional counterpart of Lemma 2.4 in Tarpey \textsl{et al.} \cite{tarpeyliflury}.

\begin{lemma}{\label{mejoraprox}}
Let $W$ be a best $k-$point approximation to $V$ and denote by $y_1, \ldots, y_k$ a set of $k$ different elements in ${\cal S}(W)$ and by ${\cal D}_j$ the domain of attraction of $y_j$. Then,
\begin{enumerate}
\item[a)] If $V\in {\cal D}_j$ then $W$ equals $y_j$ with probability 1. That is, $W = \sum_{i=1}^k y_i \uno_{{\cal D}_j} (V)$.
\item[b)]$\|V - W \| \leq \|V - y_j\|$ \textsl{a.s.} for all $y_j\in {\cal S}(W)$.
\item[c)] $E(V|W) = W$ \textsl{a.s.}, i.e., $W$ is self-consistent for $V$.
\end{enumerate}
\end{lemma}

It is worth noticing that given a self-consistent set $\{ y_1,\ldots, y_k \}$ of $V$, as in the finite--dimensional case,  we can define in a natural way a random variable $Y=\sum_{i=1}^k y_{i} \uno_{V\in {\cal D}_i}$ with support $\{ y_1, \ldots, y_k \}$ and so, $P(Y= y_j) =P(V \in {\cal D}_j)$. Since $\{ y_1, \ldots, y_k \}$ is a self-consistent set, we will have that $E(V|Y) = Y$, with probability 1. As in the finite--dimensional setting, $Y$ will not be necessarily a best approximation, unless the set $\{ y_1,\ldots, y_k \}$ is a set of $k$ principal points.


\vskip0.2in

As mentioned above, if $V$ is a random element with a self-consistent set of $k=1$ elements $\{ y_1\}$ then $E(V)=y_1$ and so that if we assume $E(V)=0$, we  have $y_1=0$. The forthcoming results  try to characterize the subspace spanned by the self--consistent points when $k>1$. They generalize the
results obtained in the finite--dimensional case by Tarpey \textsl{et al.} \cite{tarpeyliflury} and extended to gaussian processes by Tarpey and Kinateder
\cite{tarpeykini}. They also justify the use of the $k-$means algorithm not only for gaussian processes but also for elliptical processes with finite second moments.

\begin{theorem}\label{ortogoself}
Let $V$ be a random element in a separable Hilbert space $\cal H$, with finite second moment and assume that $E(V)=0$.
Let $W = \{ y_1, \ldots, y_k \}$ be a set of $k$ self--consistent points for $V$. Then,  $y_j\in Ker(\bGa_V)^{\bot}$, for all $j$, where $\bGa_V$ denotes the covariance operator of $V$.
\end{theorem}

In particular, if ${\cal W}$ denotes the linear space spanned by the $k$ self--consistent points, we have get easily that $Ker(\bGa_V) \cap {\cal W} = \{ 0 \}$. Moreover, it will also hold that $\bGa_V({\cal W}) \cap {\cal W}^{\bot} = \{ 0 \}$. This last fact will follow from the properties of semidefinite and diagonalizable operators.

\begin{corollary}\label{coro1}
Let $V$ be a random element in a separable Hilbert space $\cal H$, with finite second moment and  compact covariance operator $\bGa_V$, such that $E(V)=0$.
Let $W = \{ y_1, \ldots, y_k \}$ be a set with $k$ self--consistent points for $V$ and denote ${\cal W}$ the subspace spanned by them.
Then,
\begin{enumerate}
\item[a)] $Ker(\bGa_V) \cap {\cal W} = \{ 0 \}$.
\item[b)] Denote by  $\cal W$ be the subspace spanned by the set $\{ y_1, \ldots, y_k \}$ of $k>1$ self-consistent points. Then,
$\bGa_V{\cal W}\cap {\cal W}^{\bot}=\{0\}$.
\end{enumerate}
\end{corollary}

The following Theorems  provide the desired result relating, for elliptical elements, self--consistency and principal components.
\begin{theorem}\label{main1}
Let $V$ be a random elliptical element with $E(V)=0$ and compact covariance operator $\bGa_V$. Let $\cal W$ the subspace spanned by the set $\{ y_1, \ldots, y_k \}$ of $k>1$ self-consistent points. Then, $\cal W$ is spanned by a set of eigenfunctions of $\bGa_V$.
\end{theorem}

\begin{theorem}\label{main2}
Let $V$ be a random elliptical element with $E(V)=0$ and compact covariance operator $\bGa_V$.
If $k$ principal points of $V$ generate a subspace $\cal W$ of dimension $q$, then this subspace will also be spanned by the $q$ eigenfunctions of $\bGa_V$ related to the $q$ largest eigenvalues.
\end{theorem}

\subsection{Properties of principal points and resolution for the case $k=2$}

As mentioned above, when $k=1$  the  principal point equals the mean of the distribution. The goal of this section is  to obtain, as in the finite--dimensional setting, an explicit expression for the principal points when $k=2$. As is well known, even when dealing with finite--dimensional data, no
general result is known for any value of $k$. The following theorem will be very useful in the sequel and it generalizes a result given, for the
finite--dimensional case, by Flury \cite{flury90}. It is worth noticing that Theorem \ref{teorema52} does not require to   the random element to have an elliptical distribution.

\begin{theorem}\label{teorema52}
Let ${\cal H}$ be a separable Hilbert space and $V:\Omega\to {\cal H}$  a random element with mean $\mu$ and with $k$ principal points $\xi_1, \ldots, \xi_k \in {\cal H}$. Then, the dimension of the linear space spanned by $\xi_1 - \mu, \ldots, \xi_k - \mu$ is strictly lower than $k$.
\end{theorem}

In particular, when $k=1$ we get that  the mean is a $1-$principal point. We will now focus our attention of the case $k=2$ and the results will be derived for elliptical distributions.

\vskip0.2in
Theorem \ref{teorema35} generalizes Theorem 2 in Flury \cite{flury90} which states an analogous property for the finite dimensional vectors. As in euclidean spaces, the result assumes the existence of self-principal points for real variables, conditions under which this holds are given
in Theorem 1 of Flury \cite{flury90}.

\begin{theorem}{\label{teorema35}}
Let $V$ be an elliptical random element of a separable Hilbert space $\cal H$ with mean $\mu$ and covariance operator $\bGa$ with finite trace. Denote by $\phi_1 \in {\cal H}$ an eigenfunction of $\bGa$ with norm 1, related to its largest eigenvalue $\lambda_1$. Assume that the real random variable $Y=<v,V - \mu>$ has two principal points for any $v \in {\cal H}$ and let $\gamma_1$, $\gamma_2$ the two principal points of the real random variable $<\phi_1,V - \mu>$ . Then, $V$ has two principal points $y_1 = \mu + \gamma_1 \phi_1$ and $y_2 = \mu + \gamma_2 \phi_1$.
\end{theorem}

{\setcounter{equation}{0}
\renewcommand{\theequation}{A.\arabic{equation}}
{\setcounter{section}{0}
\renewcommand{\thesection}{\Alph{section}}

\section{Appendix}{\label{proof}}

\noi \textsc{Proof of Lemma \ref{elip:comp}}
Let $B:{\cal H}_2 \rightarrow \mathbb{R}^d$ linear and bounded, let us
show that $BAV$ is an elliptical multivariate random vector of mean
$BA\mu$ and covariance matrix $BA\bGa A^* B^*$.

Let $B\circ A: {\cal H}_1 \rightarrow \mathbb{R}^d $ the composition.
Then, $B\circ A$ is linear and bounded, therefore $BAV=(B\circ
A)(V)$ is elliptical with parameters $B\circ A(\mu)=BA\mu$ and
$(B\circ A)\bGa (B\circ A)^* = BA\bGa A^*B^*$, finishing the
proof. \square

\vskip0.2in
\noi \textsc{Proof of Lemma \ref{elip:esp}}.
a) Denote by ${\cal H}^*$ the dual space of ${\cal H}$, i.e., ${\cal H}^*$ is the set of all linear and continuous functions $f:{\cal H}\to \mathbb{R}$. Let  $f \in {\cal H}^*$, then $E(|f(V)|) <\infty$ and since $f:{\cal H} \rightarrow \mathbb{R}$ is linear and continuous it is linear and bounded. Then,    $f(V)$ has an elliptical distribution with parameters $f(\mu)$ and $f\bGa f^*$. The existence of  $E(V)$ entails that $E(f(V))$ exists and that $E(f(V))=f(E(V))$. Since $E(f(V))=f(\mu)$, by uniqueness we get that  $E(V) =\mu$.

The proof of b) will follow from the properties of the covariance operator  and Lemma \ref{elip:comp}
using the uniqueness of the covariance.

For that purpose, it will be convenient to have defined a series of special operators. Since ${\cal H}$ is separable, it  admits
an orthonormal countable base, that is, there exists $\{ \phi_n\}_{n \in \natuito}$ (eventually finite if the space is of finite
dimension) orthonormal  generating ${\cal H}$. We will choose as  basis of $\cal H$ the basis of eigenfunctions of $\bGa$ related to the eigenvalues $\lambda_1\ge \lambda_2\ge \ldots$.   {Without loss of generality we can assume that $\lambda_1>0$, otherwise $P(V=E(V))=1$ and the conclusion would be trivial.}

Define $P_n = P_{<\phi_1, \ldots, \phi_n>}$, the orthogonal projection onto the subspace ${\cal H}_1$ spanned by
$\phi_1, \ldots, \phi_n$ and $T_n$ as in (\ref{deftn}).

We want to show that $\bGa_V = {\alpha}\bGa$, i.e., that
$$<{\alpha}\bGa u, v>_{\cal H} = a_V(u,v) = Cov(<v,V>_{\cal H}, <v,V>_{\cal H})$$
for any  $ u, v \in {\cal H}$ , where we have explicitly written the space where the internal product is taken for clarity.

Let $d\in \natu$ fixed. Using that $V$ has an elliptical distribution, we get that $T_dV\sim {\cal E}_d(T_d\mu, T_d\bGa T_d^*)$. On the other hand, since  $V$ has finite second moment, the same holds for $T_dV$  which implies that  $E(T_dV)=T_d \mu$ and the covariance matrix of $T_dV$, denoted $\bSi_d$,   is {proportional} to $T_d \bGa T_d^*$. Therefore, it exists $\alpha_d\in \mathbb{R}$ such that $\bSi_d=\alpha_d T_d \bGa T_d^*$.

{We begin by showing that  $\alpha_d$ does not depend on $d$. It
is easy to see that $$T_d\bGa u=\sum_{i=1}^d \lambda_i <\phi_i,
u>_{\cal H}\;{\be_i}\;.$$ Therefore, using that
$<\phi_i, T_d^*\bx>_{\cal H}= <T_d\phi_i,
\bx>_{\realito^d}=<\be_i, \bx>_{\realito^d}=x_i$ for all
$\bx\in\mathbb{R}^d$, we obtain that
\begin{equation}
T_d \bGa T_d^*=\mbox{diag}\left(\lambda_1,\ldots,\lambda_d\right)\;.
\label{diagonal}
\end{equation}
Let $k\le d$ and $\pi_k:\mathbb{R}^d\to \mathbb{R}^k$ be the usual projection $\pi_k(\bx)=\pi_k\left((x_1,\ldots,x_p)\trasp\right)=(x_1,\ldots,x_k)\trasp=\bA_k\bx$, where $\bA_k=\left(\begin{array}{cc} \identidad_k & \bold{0}\\
\bold{0}& \bold{0}\end{array}\right)$. The fact that $\pi_k T_d V=T_k V$  implies that the covariance matrix of $T_k V$ is given by $\bSi_k=\bA_k\bSi_d\bA_k\trasp$ and so, $\alpha_k T_k \bGa T_k^*=\alpha_d \bA_k(T_d \bGa T_d^*)\bA_k\trasp$ which together with  (\ref{diagonal}), implies that $\alpha_k=\alpha_d$.

Hence, there exists $\alpha\in\mathbb{R}$ such that for all $d\in \natu$ the covariance matrix $\bSi_d$ of $T_dV$ is equal to  $\alpha T_d \bGa T_d^*$, implying that
$$<{\alpha}T_d\bGa T_d^* \bx, \by>_{\realito^d} = Cov(<\bx, T_dV>_{\realito^d} , <\by, T_dV>_{\realito^d} ).
$$
Using the definition of adjoint of $T_d$, we have that $<T_d\bGa
T_d^* \bx, \by>_{\realito^d} =<\bGa T_d^*\bx, T_d^*\by>_{\cal
H}$, meanwhile the right member of the equality can be written as
$$ Cov(<\bx, T_dV>_{\realito^d} , <\by, T_dV>_{\realito^d} )=Cov(<T_d^*\bx, V>_{\cal H}, <T_d^*\by, V>_{\cal H})\;.$$
Then, we have that for all $ d \in \natu$, $\bx$, $\by \in
\mathbb{R}^d$,
$$<{\alpha}\bGa T_d^*\bx, T_d^*\by>_{\cal H} = Cov(<T_d^*\bx, V>_{\cal H}, <T_d^*\by, V>_{\cal H})=a_V(T_d^*\bx,T_d^*\by).$$
Given $ u, v \in {\cal H}$, define $u_d=P_d u$,  $v_d=P_d v$, $\bx=T_d u=T_d u_d$ and $\by=T_d v=T_dv_d$. We have that  $\lim_{d \to \infty}\|u-u_d\|=0$ and $\lim_{d \to
\infty}\|v-v_d\|=0$. Then, using that $u_d=T_d^* \bx$, $v_d=T_d^* \by$ we get
$$ <{\alpha}\bGa u_d, v_d>_{\cal H} = <{\alpha}\bGa T_d^*\bx, T_d^* \by>_{\cal H} = a_V(T_d^* \bx, T_d^*\by) = a_V(u_d,v_d). $$
The continuity of  $a_V$ entails that $\lim_{d \to\infty}a_V(u_d,v_d) =a_V(u,v)$. On the other hand, using that $\bGa$ is a self--adjoint, compact   operator, we obtain that $\lim_{d \to \infty} <\bGa u_d, v_d>_{\cal H}=<\bGa u,  v>_{\cal H}$, which concludes the proof. \square

\vskip0.2in

\noi \textsc{Proof of Theorem \ref{eliptico_partido}.}
The proof of a) follows immediately since $\lambda_d>0$.

\noi b) It is enough to show that  $E(|f(V_2|W_1)|) < \infty$ and $E(f(V_2)|W_1) = f(\mu_2 + \bGa_{V_2, W_1}\bSi_{W_1}^{-1} (W_1 -
\mu_1))$, for any $f \in {\cal H}^*$. Let $\bW=(W_1,f(V_2))=TV$. Using that with $T$ is a bounded and linear operator, we get that
$\bW$ is elliptical of parameters $(\bmu_1, f(\mu_2))\trasp$ and
$$T\bGa T^*=\left(
\begin{array}{cc} \bSi_{W_1} & Cov(W_1,f(V_2))\\
Cov(W_1,f(V_2))\trasp & f\bGa f^*
\end{array}\right)\;.$$
Using Theorem \ref{elip:finito},  we get that  $f(V_2)|W_1$  has also an elliptical distribution with expectation given by $\bmu_f=f(\mu_2) + Cov(f(V_2),W_1)\bSi_{W_1}^{-1} (W_1 - \bmu_1)$.
On the other hand, $Cov(f(V_2),W_1)=f\bGa_{V_2, W_1}$ which implies that
\begin{eqnarray*}
E(f(V_2)|W_1)&=&f(\mu_2) + Cov(f(V_2),W_1)\bSi_{W_1}^{-1} (W_1 - \bmu_1)\\
&=&f(\mu_2) + f\bGa_{V_2, W_1}\bSi_{W_1}^{-1} (W_1 -
\bmu_1)=f(\mu_2 + \bGa_{V_2,W_1}\bSi_{W_1}^{-1} (W_1 - \bmu_1))
\end{eqnarray*}
and so, we conclude the proof. \square

\vskip0.2in
\noi \textsc{Proof of Lemma \ref{selfunitario}}.
a) Using that $\cal W$ is self-consistent for $V$, we get that $E(V | V \in {\cal D}_j) = y_j$. Let us notice that, since $U$ is a unitary operator,  $V \in {\cal D}_j$ if and only if $V_2 = \nu + \rho\; UV\in \widetilde{{\cal D}_j}$. Therefore, $\widetilde{{\cal
D}_j}$ is the domain of attraction of $\nu + \rho\; U y_j$ which implies that  $\nu + \rho\; U{\cal D}_j\subset \widetilde{{\cal
D}_j} $. Hence,
\begin{eqnarray*}
E(V_2 | V_2 \in \widetilde{{\cal D}_j}) &=& E(\nu + \rho UV | V_2 \in \widetilde{D_j}) = \nu + \rho UE(V | \nu + \rho UV \in \nu + \rho U{\cal D}_j)\\
&=& \nu + \rho UE(V | V \in {\cal D}_j) = \nu + \rho Uy_j.
\end{eqnarray*}
b) Let $\xi_1, \ldots, \xi_k$ be any set of points in $\cal H$ and denote by ${\cal A}_1, \ldots, {\cal A}_k$ their respective domains of attraction. We have to prove that $E(d^2(V_2, {\cal W}_2 )) \leq E(d^2(V_2,\{\xi_1, \ldots, \xi_k\}))$.

Let $z_j$ such that $\xi_j = \nu + \rho \; U z_j$, then
\begin{eqnarray*}
E(d^2(V_2,\{\xi_1, \ldots, \xi_k\})) &=& E( \min_{1\leq j \leq k} \|V_2 - \xi_j\|^2)= E( \min_{1\leq j \leq k} \|\nu + \rho\; UV - \xi_j\|^2) \\
&=&E( \min_{1\leq j \leq k} \|\nu + \rho \;UV - \nu - \rho\; U z_j\|^2) =E( \min_{1\leq j \leq k} \|\rho \;UV - \rho\; U z_j\|^2)\\
&=&\rho^2 E( \min_{1\leq j \leq k} \|UV - U z_j\|^2)=\rho^2 E(
\min_{1\leq j \leq k} \|V -  z_j\|^2)
\end{eqnarray*}
where the last inequality holds from the fact that $U$ is an isometry. On the other hand, using that $\cal W$ is a set of principal points of $V$, we get that  $\{y_1,\ldots,y_k\}=\argmin_{z_1,\ldots ,z_k} E( \min_{1\leq j \leq k} \|V -  z_j\|^2)$. Threfore, we have that
\begin{eqnarray*}
E(d^2(V_2,\{\xi_1, \ldots, \xi_k\}))&=&\rho^2  E( \min_{1\leq j \leq
k} \|V -  z_j\|^2)\ge \rho^2  E( \min_{1\leq j \leq k} \|V -
y_j\|^2)=E( \min_{1\leq j \leq k} \|V_2 - \xi_{0,j}\|^2)
\end{eqnarray*}
where $\xi_{0,j} = \nu + \rho Uy_j$, which means ${\cal W}_2$ are principal points of $V_2$. Besides, we also obtain that $E(d^2(V_2,{\cal W}_2)) = \rho^2 E(d^2(V,{\cal W}))$. \square

\vskip0.3in

\noi \textsc{Proof of Lemma \ref{finitedim}}.
Let us define $A:{\cal H} \rightarrow \mathbb{R}^q$ as $A(v) =\bv=(v_1,\ldots,v_q)\trasp$ with $v_i= <e_i, v>$.
We want to show that $\bX=AV \in \mathbb{R}^q$ has as a self-consistent set
${\cal W}=\{\bw_{1}, \ldots, \bw_k\}=\{ Ay_1,
\ldots, Ay_k \}$. Denote by $\widetilde{{\cal D}_j}$ the
domain of attraction of $Ay_j$ and by ${\cal D}_j$ that of $y_j$. Then, by extending the orthonormal basis $\{e_1, \ldots, e_q\}$ of $\cal M$ to an orthonormal basis $\{e_{\ell}, \ell\ge 1\}$ of $\cal H$, and using that $<y_j,e_{\ell}>=0$ for $\ell>q$, it is easy to see that  $V\in {\cal D}_j$ if and only if $AV\in\widetilde{{\cal D}_j}$, which implies that
$$E(AV | AV \in \widetilde{{\cal D}_j}) = A E(V | AV \in \widetilde{{\cal D}_j}) = A E(V | V \in {\cal D}_j) = Ay_j=\bw_j\;, $$
concluding the proof. \square

\vskip0.2in

\noi \textsc{Proof of Lemma \ref{mejoraprox}}. a) We will always suppose that the probability of $V$ being found in the frontier of two domains of attractions (${\cal D}_i \cap {\cal D}_j$) is 0.

Let us suppose that the result is false and denote by $\Omega$ the common probability space. That is, let us assume that the set $A=\{\omega: V(\omega) \in {\cal D}_j \} \cap \{ W(\omega) = y_r \} $ has probability strictly positive, with $r\neq j$ and define a new element $Z(\omega)$ that is equal to $W(\omega)$ if $\omega\notin A$ and is equal to $y_i$ if $\omega\in A$. We will show that $Z$ is a better approximation of $V$ than $W$.
\begin{eqnarray*}
E(\|V-Z\|^2) &=& E( \|V-Z\|^2 \uno_{A^c} ) + E( \|V-Z\|^2\uno_{A}  ) \\
&=& E( \|V-W\|^2 \uno_{A^c} ) + E( \|V-y_j\|^2 \uno_{A}  ) \\
\end{eqnarray*}
However, since for any $\omega\in A$, we have that $V(\omega) \in {\cal D}_j$, we get that $
\|V-y_j\|^2 \uno_{A} < \|V-y_r\|^2 \uno_{A} = \|V-W\|^2 \uno_{A}
$. Therefore, using that $A$ is a set with  positive probability,  we obtain that $ E( \|V-y_j\|^2 \uno_{A}  ) < E( \|V-W\|^2 \uno_{A}  )$ and so
\begin{eqnarray*}
E( \|V-W\|^2 \uno_{A^c} ) + E( \|V-y_j\|^2 \uno_{A}) < E( \|V-W\|^2 \uno_{A^c} ) + E( \|V-W\|^2 \uno_{A}) =
E( \|V-W\|^2)\; ,
\end{eqnarray*}
which entails that  $ E(\|V-Z\|^2) < E( \|V-W\|^2)$ implying that $Z$ is a better $k-$point approximation than $W$, concluding the proof of a).

\noi b) If $V(\omega) \in {\cal D}_j $ then, using a) we get that, except for a zero probability set, $W(\omega) = y_j$. Then, $\|V -W\| (\omega) = \|V(\omega)-y_j\|\le\|V(\omega)-y_i\|$ for any $y_i$ since $V(\omega) \in {\cal D}_j$.

\noi c)  $E(V|W)$ is  a measurable function $g(W)$ that minimizes
the expected squared distance between $V$ and any measurable
function $h(W)$. Using a) we have that, any function $h(W)$  has
a  support containing at most $k$ points. Then, by the definition
of best approximation, we  have that {$W$ is a
better approximation than $h(W)$ for any measurable function $h$},
with the expected squared distance criteria. Then, the function
$g(W)=E(V|W)$ equals $W$. \square

\vskip0.2in
In order to prove Theorem \ref{main1} we will need some technical Lemmas. In particular, these lemmas will allow to derive that the matrix
$\bGa_{W_1, W_1}$ defined therein is not singular, which is a fundamental step in order to get the desired conclusion.

\begin{lemma}\label{lemaapendice}
Let $V$ be a random element in a separable Hilbert space $\cal H$, with finite second moment and assume that $E(V)=0$.
\begin{enumerate}
\item[a)] If $h \in Ker(\bGa_V)$ then $<h,V>\equiv 0 \in \mathbb{R}$ \textsl{a.s.}
\item[b)] Denote by ${\cal H}_1 = Ker(\bGa_V)$, the kernel of the covariance operator and ${\cal H}_2 = Ker(\bGa_V)^{\bot}$ its orthogonal. Then,  $P\left(P_{{\cal H}_1}V = 0\right)=1$   and  $P\left(P_{{\cal H}_2} V = V\right)=1$, i.e, $P\left(V \in Ker(\bGa_V)^{\bot} \right)= 1$.
\end{enumerate}
\end{lemma}
\noi \textsc{Proof}. The proof of a) follows easily noticing that $E(<h,V>)= <h,E(V)> = 0$ and  $Var(<h,V>) = Cov(<h,V>, <h,V>) = <\bGa_V h, h> = 0$, since $h \in Ker(\bGa_V)$.

b) Note that the separability of $\cal H$ entails that $ V = P_{{\cal H}_1} V + P_{{\cal H}_2} V$, since $ Ker(\bGa_V)$ is a closed subspace. Thus, it will be enough to show that $P_{{\cal H}_1} V=0$ with probability $1$. For the sake of simplicity, we will assume that both ${\cal H}_1$ and ${\cal H}_2$ are  infinite dimensional spaces. Otherwise, the same calculations hold  but using a finite index set as the only significant change. Denote by $\{e_1, e_2, \dots\}$ an orthonormal basis of ${\cal H}_1$ and extend it to a basis of $\cal H$ so that, $\{f_1, f_2, \dots\}$ will denote  the orthonormal basis of  $ {\cal H}_2 $. Then, $Var(P_{{\cal H}_1} V) = E\left(\|P_{{\cal H}_1} V-E\left(P_{{\cal H}_1} V\right)\|^2\right)=\sum_{i=1}^{\infty} <\bGa_{P_{{\cal H}_1} V} \; e_i, e_i>=\sum_{i=1}^{\infty} <\bGa_{V}\; e_i, e_i>$ and each summand equals zero since $e_i \in Ker(\bGa_V)$. On the other hand, $E(P_{{\cal H}_1} V) = P_{{\cal H}_1} E(V) = 0$, therefore  $E\left(\|P_{{\cal H}_1} V\|^2\right)=0$, which implies that $P_{{\cal H}_1} V=0$ \textsl{a.s.}, concluding  the proof. \square

{\begin{corollary}
If $A\cap Ker(\bGa_V)^{\bot} = \emptyset$ then $P(V\in A) = 0$.
\end{corollary}}

The proof is immediately since  $P(V \in Ker(\bGa_V)^{\bot} = 1$.

\vskip0.2in
\noi \textsc{Proof of Theorem \ref{ortogoself}}. Note that $Ker(\bGa_V)^{\bot}$ is a closed subspace and therefore a convex set. For each point $y_j$, its domain of attraction  ${\cal D}_j$  is also a convex set, therefore, $Ker(\bGa_V)^{\bot}\cap {\cal D}_j$  is also convex.
By Lemma \ref{lemaapendice},  the support of the random element $V$ is included in  $Ker(\bGa_V)^{\bot}$, thus
$$
y_j = E(V| V \in {\cal D}_j) = E(V| V \in {\cal D}_j, V \in
Ker(\bGa_V)^{\bot}) = E(V | V \in {\cal D}_j \cap
Ker(\bGa_V)^{\bot} ).
$$
Now the proof follows easily by noticing that the expectation of a random element taking values in  a convex set $\cal C$ will also be in $\cal C$, i.e., $y_j \in {\cal D}_j \cap Ker(\bGa_V)^{\bot} \subset Ker(\bGa_V)^{\bot}$.  \square

\vskip0.2in

\noi \textsc{Proof of Corollary \ref{coro1}}. a) follows immediately from Theorem \ref{ortogoself}. We have only to prove b).
Let $z \in \bGa_V({\cal W}) \cap {\cal W}^{\bot}$, then $z = \bGa_V w$ with $w \in {\cal W}$. We want to show that $z=0$. Let $\{ \phi_j
\}_{j \in \mathbb{N}}$ be an orthonormal base of eigenfunctions of $\bGa_V$ related to the eigenvalues $\mu_1\ge  \ldots\ge \mu_j\ge \ldots$
Then, $ w = \sum_{j=1}^{\infty} <w, \phi_j>\phi_j$ which entails that $z =  \sum_{j=1}^{\infty} <w, \phi_j>\mu_j \phi_j$. Using that $z \in {\cal W}^{\bot}$, we get that  $<z, w> = 0$ and so
$ 0 = <z,w> = \sum_{j=1}^{\infty} \mu_j <w, \phi_j>^2$. The fact that $\mu_j$ are non-negative, implies that $<w, \phi_j>=0$ if $\mu_j>0$ and so $\bGa_V w = 0$, which entails that $w \in Ker(\bGa_V)$ and $z= 0$ concluding the proof. \square

\vskip0.2in

\noi \textsc{Proof of Theorem \ref{main1}}. For the sake of simplicity, we will avoid the index $V$ in $\bGa_V$ and will denote $\bGa$ the covariance operator.

Let $q$ be the dimension of ${\cal W}$ and   $\{ v_1, \ldots, v_q \}$ an orthonormal basis of $\cal W$. Let us denote by  $\{ v_1, \ldots, v_q, v_{q+1}, \dots, \}$ the extension to an orthonormal  basis of $\cal H$. Define $A_1: {\cal H} \rightarrow
\mathbb{R}^q$ as $A_1(h) = \sum_{i=1}^q <v_i,h>\be_i$,
with $\be_i$ the canonical basis of $\mathbb{R}^q$. Then,
$A_1^{*}:\mathbb{R}^q\to {\cal H}$, equals $A_1^*\bx=\sum_{i=1}^q x_i v_i$ with $\bx=(x_1,\ldots,x_q)\trasp$, that is, the image of
$A_1^*$ is $\cal W$. We also define $A_2:
{\cal H} \rightarrow {\cal H}_{\infty}\subset \mathbb{R}^{\mathbb{N}}$ as
$$
A_2 (h) = \left( \begin{array}{c}
<v_{q+1}, h>   \\
<v_{q+2}, h>   \\
\vdots
\end{array} \right).
$$
To ensure the continuity of the second operator, we consider as norm in ${\cal H}_{\infty}$  the norm given by the square root of the sum of squares of the elements of the sequence and as inner product, the one generating this norm. Using  Parseval's identity we have that for any $h \in {\cal H}$, $\|h\|^2 =\sum_{k=1}^{\infty} <h, v_k>^2$ which implies that $A_2$ is continuous and with norm equal to 1 since
$$\|A_2(h)\|^2 = \sum_{k=q+1}^{\infty}<h,v_k>^2 \leq \sum_{k=1}^{\infty} <h, v_k>^2 = \|h\|^2\;,$$
with equality for any $h \in {\cal W}^{\bot}$. Moreover, define $A: {\cal H} \rightarrow {\cal H}_{\infty}$, as $A(h)=\left(A_1(h)\trasp, A_2(h)\trasp\right)\trasp$ and $W=A(V)=\left(A_1(V)\trasp, A_2(V)\trasp\right)\trasp=\left(W_1\trasp,W_2\trasp\right)\trasp$ with $W_1$ finite--dimensional. Notice that $W$ is an elliptical element in ${\cal H}_{\infty}$ with null expectation  and covariance operator $\bGa_{AV} = A\bGa_{V} A^*$.  Thus, if $\bGa_{W_1, W_1}$ is non--singular, using Theorem \ref{eliptico_partido}, we get that
\begin{equation}
E( A_2(V)| A_1(V)) = E(W_2|W_1 ) = \bGa_{W_2,W_1}(\bGa_{W_1,W_1})^{-1}A_1(V)
\label{espcond}
\end{equation}
where $\bGa_{W_2,W_1} = A_2\bGa A_1^*$ and $\bGa_{W_1,W_1} =A_1\bGa A_1^*$.

Using that  $Ker(\bGa_V)\cap {\cal W}=\{0\}$ and that $\bGa_V{\cal W}\cap {\cal W}^{\bot}=\{0\}$, we will show that $\bGa_{W_1, W_1}$ is non--singular.
Since  $\bGa_{W_1, W_1}$ is an endomorphism between finite dimensional vector spaces, we only have to prove that it is a monomorphism
and we will automatically have that it is an isomorphism. Let us see the injectivity of this operator. Let us assume that for some $\bh \in
\mathbb{R}^q$, we have that $A_1\bGa A_1^*\bh = \bf{0}$. We want to show that $\bh = \bf{0}$. Using that $A_1\bGa A_1^*\bh = \bf{0}$, we get that $\bGa A_1^*\bh \in Ker(A_1) = {\cal W}^{\bot}$. On the other hand, $A_1^*\bh \in {\cal W}$ and so,  $\bGa A_1^*\bh \in\ \bGa({\cal W}) \cap {\cal W}^{\bot}$ which implies $\bGa A_1^*\bh = 0$ by Corollary \ref{coro1}. Hence, $A_1^*\bh \in Ker(\bGa) \cap {\cal W}=\{ 0 \}$ by Corollary \ref{coro1} and so, $A_1^* \bh = 0$. The fact that $A_1^*$ is injective, leads to  $\bh = \bf{0}$ and therefore, $\bGa_{W_1, W_1}$ is non--singular and (\ref{espcond}) holds.

Define now  a random element $Y$ such that $P(Y=y_j) = P(V\in {\cal D}_j)$, with ${\cal D}_j$ the domain of attraction of $y_j$, that is, $Y=\sum_{i=1}^k y_{i} \uno_{V\in {\cal D}_i}$. Using that $A_2$ is linear and continuous, and that $y_j$ is an element on the self--consistent set, we get that
$$E(A_2 V | Y = y_j) = E(A_2 V | V \in {\cal D}_j) = A_2E(V | V \in {\cal D}_j) = A_2y_j=0\;$$
where the last equality holds since $y_j\in {\cal W}$, that is the kernel of $A_2$. Therefore,  $P(E(A_2V | Y)=0)=1$ and so, $0=E(A_2V| Y)=E(E(A_2V|A_1V)|Y)$ with probability 1. Then, using (\ref{espcond}), we get
\begin{eqnarray*}
0 = E(A_2V| Y) &=&E( \bGa_{W_2, W_1}(\bGa_{W_1,W_1})^{-1}A_1(V) | Y) = \bGa_{W_2,W_1} (\bGa_{W_1,W_1})^{-1} A_1\;E(V|Y)\\
&=& A_2\bGa A_1^* (A_1\bGa A_1^*)^{-1}A_1\;Y \quad \mbox{\textsl{a.s.}}
\end{eqnarray*}
where the last equality follows using Lemma \ref {mejoraprox}.

Using that the support of $Y$ spans $\cal W$, we get that the
support of $A_1Y$ spans $\mathbb{R}^q$.  Then, using that $ (A_1\bGa
A_1^*)^{-1}$ is non--singular, the fact that $ P\left(A_2\bGa A_1^* (A_1\bGa
A_1^*)^{-1}A_1 Y=0\right)=1$ implies that $A_2\bGa A_1^* (A_1\bGa
A_1^*)^{-1}A_1\; y_j=0 $, for $1\le j\le k$, and so, $\bGa{\cal
W}\cap {\cal W}^\bot=\{0\}$. Therefore, $\forall \bx\in
\mathbb{R}^q$, $A_2\bGa A_1^*\bx=0$, or equivalently, $A_2\bGa
A_1^*: \mathbb{R}^q \rightarrow {\cal H}_{\infty}$ is the null
operator and the same will be   true for $A_1\bGa A_2^*:{\cal
H}_{\infty}\to \mathbb{R}^q$.

Define the projection operators $P_{\cal W} = A_1^* A_1:{\cal H}\to {\cal W}$ (the projection over ${\cal W}$) and $P_{{\cal W}^\bot} = A_2^*A_2:{\cal H}\to {\cal W}^\bot$ (the projection over ${\cal W}^\bot$). Then,  $P_{{\cal W}^\bot}\bGa P_{{\cal W}}=0$ and $P_{{\cal W}}\bGa P_{{\cal W}^\bot}=0$ and so,  $\bGa {\cal W} \subset {\cal W}$ and  $\bGa  {\cal W}^\bot) \subset {\cal W}^\bot $ which  implies that ${\cal W}$ and ${\cal W}^\bot$ are $\bGa-$invariant, i.e., ${\cal W}$ decomposes $\bGa$.

Then, the restriction of the covariance operator  $\bGa|_{{\cal W}}: {\cal W} \rightarrow {\cal W}$  will be well defined. $\bGa$ is compact and self-adjoint and so, diagonalizable. Besides, $\bGa$ restricted to ${\cal W}$ will also be compact and self-adjoint and thus, it will  be diagonalizable with the same eigenfunctions. We will then have that ${\cal W}$, the domain of $\bGa|_{\cal W}$ is spanned by a set of eigenfunctions of $\bGa$. Using that ${\cal W}$ has dimension $q$, we get that ${\cal W}$ is spanned by $q$ eigenfunctions of $\bGa$ concluding the proof. \square

\vskip0.2in

\noi \textsc{Proof of Theorem \ref{main2}}. As in the {proof of Theorem \ref{main1}}, we will avoid the index $V$ in $\bGa_V$. Let $\lambda_1 \ge\lambda_2 \ldots $ the ordered eigenvalues of $\bGa$, with its corresponding eigenfunctions $\phi_1,\phi_2, \ldots $. Let $\{ y_1, \ldots, y_k \}$ be a set of $k$ principal points that span $\cal W$. Theorem \ref{main1} entails that $\cal W$ is spanned by $q$ eigenfunctions of $\bGa$. Let $r$ be an integer such that
$\{ \phi_1, \ldots, \phi_r \}$ contains the $q$ eigenfunctions that generate $\cal W$.

Denote, for each principal point $y_j$, $a_{ji}=<\phi_i,y_j>$, so that, $a_{ji}=0$ for $j>r$ and  $y_j=\sum_{i=1}^r a_{ji}\phi_i$.
Define $\bold{a}_j = (a_{j1}, \ldots, a_{jr})\trasp$ and $\bX
=(<\phi_1, V>, \ldots, <\phi_r, V>)\trasp.$

The eigenfunctions will conform an orthonormal basis, then, we  have that
\begin{eqnarray*}
\|V - y_j\|^2 &= & \| \sum_{i=1}^{r}(<\phi_i,
V>-a_{ji})\phi_i\|^2 + \| \sum_{i=r+1}^{\infty}<\phi_i,
V>\phi_i \|^2
\\
&=& \sum_{i=1}^r(<\phi_i,V> - a_{ji})^2 + \sum_{i=r+1}^{\infty}
<\phi_i,V>^2 = \|\bX - \ba_j\|_r^2 + \sum_{i=r+1}^{\infty} <\phi_i,V>^2\; ,
\end{eqnarray*}
where $\| \cdot \|_r$ is the euclidean norm in $\mathbb{R}^r$. Note that the $k$ principal points are the $k$ points that minimize $MSE(V, \{ \xi_1, \ldots, \xi_k \})=E(\min_{1\le j\le k} \|V - \xi_j\|^2)$ over the sets of $k$ points. Besides,
$$MSE(V, \{ y_1, \ldots, y_k \}) = E(\min_{1\le j\le k} \|V - y_j\|^2) =E(\min_{1\le j\le k}\|\bX - \ba_j\|_r^2) + \sum_{i=r+1}^{\infty}\lambda_i.
$$
and so, $MSE(\bX, \{ \ba_1, \ldots, \ba_k \}) = MSE(V, \{ y_1, \ldots, y_k\}) - \sum_{i=r+1}^{\infty}\lambda_i$.

Using that $\{ y_1, \ldots, y_k \}$  minimize $MSE(V, \cdot)$, it is easy to obtain that $\{ \ba_1, \ldots, \ba_k \}$ minimize $MSE(\bX, \{
\bb_1, \ldots, \bb_k \})$, over the sets of $k$ points in $\mathbb{R}^r$ which entails that $\{ \ba_1, \ldots, \ba_k \})$ is a set of $k$
principal points for $\bX$. On the other hand,  $\bX$ has an elliptical distribution since  $V$ has an elliptical distribution. Thus, using the result in Tarpey, Li and Flury  \cite{tarpeyliflury}, we obtain that the $k$ principal points of $\bX$ lie in the linear space related to the $q$ largest eigenvalues of the covariance, $\bSi_{\bX}$, of $\bX$.

Define $A:{\cal H}\rightarrow \mathbb{R}^r$ as $A(v) = (<\phi_1, v>, \ldots, <\phi_r,v>)\trasp$. $A$ is a linear and bounded operator. Denote by $A^{*}: \mathbb{R}^r \rightarrow {\cal H}$ the adjoint operator, i.e., $ A^{*}(\bx) = \sum_{i=1}^r x_i \phi_i$. Noticing that $\bX = AV$ and that $\ba_j = Ay_j$, we get easily that $\bSi_{\bX} =A\bGa A^{*}$. Therefore, the $q$ largest eigenvalues of $\bSi_{\bX}$ will be  equal to the $q$ largest eigenvalues of
$\bGa$. Moreover, the eigenvectors of $\bSi_{\bX}$ can be written as $A\phi_i$.

In conclusion, $\{ \ba_1, \ldots, \ba_k \} =\{ Ay_1, \ldots, Ay_k \}$ are the $k$ principal points of $\bX =AV$ and the linear space spanned by them is spanned by $A\phi_i$ with $i=1, \ldots, q$. By restricting $A$ to the space linear space ${\cal M}$ spanned by $\{\phi_1,\ldots, \phi_r\}$ we  have a surjective isometry and so,   if we define $\tilde{A}: {\cal M} \rightarrow \mathbb{R}^r$ as the restriction of $A$ to the subspace ${\cal M}$,
its inverse will be given by $\tilde{A^{*}}: \mathbb{R}^r \rightarrow{\cal M}$ which is essentially equal to $A^{*}$ except for the codomain. The proof follows now easily by noticing that  since $y_j\in {\cal M}$,   $\{ Ay_1,
\ldots, Ay_k \} = \{ \tilde{A}y_1, \ldots, \tilde{A}y_k \}\subset \tilde{\cal M}$ with $\tilde{\cal M}$ the linear space spanned by $\{\tilde{A}\phi_1, \ldots,\tilde{A}\phi_q\}$. Hence, applying $\tilde{A}^{-1} = \tilde{A^{*}}$, we  get that $\{ y_1, \ldots, y_k \}$ is included in the linear space spanned by $\{\phi_1, \ldots, \phi_q\}$. \square

\vskip0.2in

\textsc{Proof of Theorem \ref{teorema52}}. Without loss of generality we can assume that $V$ has mean $\mu=0$.
Let $c_1, \ldots, c_k \in {\cal H}$ be arbitrary elements and define $b_i = c_i - c_k$, $1\le i\le k-1$. Denote by $m\leq k-1$ the dimension of the linear space ${\cal M}_1$ spanned by $b_1,\ldots, b_{k-1}$ and by ${\cal M}_2 = {\cal M}_1^{\bot}$. Let  $a_1, \ldots, a_m$ be an orthonormal basis of ${\cal M}_1$ so that,  ${\cal M}_2$ will be spanned by $a_{m+1}, a_{m+2}, \ldots $, being $\{a_i\}_{i\ge 1}$ an orthonormal basis of $\cal H$.

As in (\ref{deftn}), define  $A_1:{\cal H} \rightarrow {\cal H}_{\infty}\subset\mathbb{R}^{\mathbb{N}}$ as
$$
A_1 (a_i) = \left\{ \begin{array}{cc}
(e_i)_{\mathbb{R}^{\mathbb{N}}} & i =1,2, \ldots , m  \\
0 & i > m
\end{array} \right.
$$
and $A_2:{\cal H} \rightarrow {\cal H}_{\infty}\subset \mathbb{R}^{\mathbb{N}}$ as
$$
A_2 (a_i) = \left\{ \begin{array}{cc}
0 & i =1,2, \ldots , m  \\
(e_i)_{\mathbb{R}^{\mathbb{N}}} & i > m
\end{array} \right.
$$
Furthermore, let $A:{\cal H} \rightarrow {\cal H}_{\infty}\subset\mathbb{R}^{\mathbb{N}}$ be $A(h) = A_1(h) + A_2(h)$. Notice that $A(a_i) = e_i$ so that $A$ is a surjective isometry and so, since it is an unitary application, its inverse will be its adjoint, the $A^*:{\cal H}_{\infty} \rightarrow {\cal H}$ such that $A^*(e_i) = a_i$. Define $d_i = Ac_i = A_1c_i + A_2c_i = d_i^{(1)} + d_i^{(2)}$. Then, using that $c_i - c_j=b_i-b_j\in {\cal M}_1$, we get that $d_i^{(2)} - d_j^{(2)}= A_2c_i - A_2c_j = A_2(c_i - c_j)=0$, that is, all the values $d_i^{(2)}$, $1\le i \le k$, are equal to a value that we will denote by $d^{(2)}$. Moreover, we have that $d_i = d_i^{(1)} + d^{(2)}$, with both  terms orthogonal between themselves.

Define $W_1 = A_1V$, $ W_2 = A_2V$, $W = W_1 + W_2 = A V$ and $W_1$ and $W_2$ are orthogonal. Using that  $A$ is unitary, we get that
\begin{equation}
E(d^2(V, \{c_1, \ldots, c_k\})) = E(d^2(AV, \{Ac_1, \ldots, Ac_k\})) =
E(d^2(W, \{d_1, \ldots, d_k\})). \label{ed}
\end{equation}
Then, $\|W - d_i\|^2 = \| W_1 - d_i^{(1)}\|^2 + \|
W_2 - d^{(2)}\|^2$,
and so that (\ref{ed}) can be written as
$$
E(d^2(V,\{c_1, \dots, c_k\}))= E(d^2(W_1, \{ d_1^{(1)}, \dots,d_k^{(1)}\})) + E(d^2(W_2,\{ d^{(2)}, \dots, d^{(2)}\})).
$$
The second term on the right hand side  equals $E(d^2(W_2, \{d^{(2)}\}))$ which is minimized when $d^{(2)}= E(W_2) = A_2\mu = 0$. Therefore,
$$
E(d^2(V,\{c_1, \ldots, c_k\}))\geq E(d^2(W_1,\{ d_1^{(1)}, \dots,d_k^{(1)}\})) + E(d^2(W_2,\{ 0, \dots, 0\})),
$$
reaching the equality when $d^{(2)} = 0$. Define $c_i^* = A_1^*d_i^{(1)} =A_1^*A_1c_i =P_{{\cal M}_1}c_i $, $V_1 = A_1^*W_1$  and $V_2 = A_2^*W_2$ ($V = V_1 + V_2$), then
\begin{eqnarray*}
E(d^2(V,\{c_1, \ldots, c_k\})) &\geq& E(d^2(W_1, \{d_1^{(1)}, \ldots, d_k^{(1)}\})) + E(d^2(W_2, \{0, \ldots, 0\}))\\
&=& E(d^2(A_1^*W_1, \{A_1^*d_1^{(1)}, \ldots, A_1^*d_k^{(1)}\})) + E(d^2(A_2^*W_1, \{0, \ldots, 0\}))\\
&=& E(d^2(V_1, \{A_1^*d_1^{(1)}, \ldots, A_1^*d_k^{(1)}\})) + E(d^2(V_2, \{0, \ldots, 0\}))\\
&=& E(d^2(V_1, \{A_1^*A_1c_1, \ldots, A_1^*A_1c_k\})) + E(d^2(V_2,\{ 0,\ldots, 0\}))\\
&=&E(d^2(V,\{A_1^*A_1c_1 + 0, \ldots, A_1^*A_1c_k + 0\})) = E(d^2(V,\{c_1^*, \ldots, c_k^*\})) \;.
\end{eqnarray*}
where the last equality follows from the orthogonality of the decomposition. Summarizing \linebreak $E(d^2(V,\{c_1^*, \ldots, c_k^*\})) \leq E(d^2(V,\{c_1, \ldots, c_k\}))$, where the equality holds if $c_i = A_1^*A_1c_i = P_{{\cal M}_1}c_i$.
At the principal points we will get the equality since by definition principal points  minimize $E(d^2(V,\{c_1, \ldots, c_k\}))$, hence,  if   $c_i$
correspond to the principal points  $\xi_i$,  then $c_i\in{\cal M}_1$. Using that  ${\cal M}_1$ has dimension lower or equal than $k-1$, we obtain the desired result. \square

\vskip0.2in

The following result, which we state for completeness, can be found in Flury \cite{flury90}.

\begin{lemma} \label{lema:aux:va} Let $Y_1$ and $Y_2$ be two real random variables such that $Y_2$ has the same distribution as $\rho Y_1$
for some value of $\rho$. Then,
$$
D_{Y_1}(k) / Var(Y_1) = D_{Y_2}(k) / Var(Y_2)\;.
$$
\end{lemma}


\vskip0.2in
\textsc{Proof of Theorem \ref{teorema35}} Without loss of generality, we sill assume that $\mu=0$. So as to reduce notation burden, define
$D_V(c_1, c_2) = E(d^2(V,\{c_1, c_2\}))$.

We will first show that $D_V$ is minimized if the two elements $c_1, c_2 \in {\cal H}$  lie on a straight line with direction $c_2 - c_1$. Theorem
\ref{teorema52} allows us to do this. Effectively, in the proof of Theorem \ref{teorema52}, we derived that each   principal point (assuming existence) $y_i$ belongs to the linear space ${\cal M}_1$  spanned by $y_2 - y_1$. That is, both elements lie in a straight line with direction $a_1=(y_2 - y_1)/\|y_2 - y_1\|$.

Take $c_1, c_2 \in {\cal H}$, $c_1\ne c_2$ and let ${\cal M}_1$ the linear space of dimension $1$ spanned by $a_1=(c_2 - c_1)/\|c_2 - c_1\|$ and ${\cal M}_2={\cal M}_1^\bot$ with orthonormal base $\{a_j: j\ge 2\}$. So, using the same notation as in Theorem \ref{teorema52}, we consider $A_1:{\cal H} \rightarrow {\cal H}_{\infty}$ defined as $A_1 (a_j)=0$ if $j\ge 2$ and $A_1 (a_1) =e_1$ with  $e_j$ the element of ${\cal H}_{\infty}$ with its $j$ coordinate equal to $1$ and all the others equal to 0, and $A_2:{\cal H} \rightarrow {\cal H}_{\infty}$
defined as $A_2 (a_1) =0$, $A_2(a_i)=e_i$ if $i > 1$. Let $W_1 = A_1V=<a_1,V> e_1$, $W_2 = A_2V$, $W=W_1+W_2$. Using that ${\cal M}_1$ is a one dimensional subspace, we get that $W_1$  has the same distribution as the random variable $Y_1=<a_1,V>$, which is elliptic and so symmetric around $0$, since $V$ is elliptic.  Let us remember that if $d_i=Ac_i=A_1c_i+A_2c_i=d_i^{(1)}+d_i^{(2)}$, then $d_1^{(2)}=d_2^{(2)}=d^{(2)}$ and
\begin{eqnarray*}
E(d^2(V,\{c_1, c_2\})) &=& E(d^2(W, \{d_1, d_2\})) = E(d^2(W_1,\{ d_1^{(1)},d_2^{(1)}\})) + E(d^2(W_2, \{d^{(2)},d^{(2)}\}))\\
&\ge  & E(d^2(W_1, \{d_1^{(1)},d_2^{(1)}\})) + E(d^2(W_2, \{0,0\})).
\end{eqnarray*}
Then, for fixed $c_1$ and $c_2$, $E(d^2(W_1, \{d_1^{(1)}, d_2^{(1)}\}))$ can be minimized taking $d_1^{(1)}$ and $d_2^{(1)}$
as the principal points of $Y_1=<a_1,V>$. Let $\xi_1$ and $\xi_2$ be the principal points of $Y_1$.
Define $d_1^* = \xi_1 e_1$, $d_2^* = \xi_2 e_1$. It follows that $E(d^2(W,\{d_1^*, d_2^*\})) \leq E(d^2(W,\{d_1, d_2\}))$,
with equality if $d_1 = d_1^*$ and $d_2 = d_2^*$.

Using that $W = AV$, we get $V = A^*W$ and so,
$$
c_1^* = A^*d_1^* = A^* \xi_1 e_1 = \xi_1 A_1^* e_1=\xi_1 a_1= \xi_1 (c_2 - c_1) / \|c_2 - c_1\|.
$$
Analogously,  $c_2^* = \xi_2 (c_2 - c_1) / \|c_2 - c_1\|$. Hence,
$$
E(d^2(V,\{\xi_1 a_1, \xi_2 a_1\})) = E(d^2(V,\{c_1^*, c_2^*\})) \leq
E(d^2(V,\{c_1, c_2\}))\;.
$$
Given $a \in {\cal H}$ such that $\|a\| = 1$, for each pair $c_1, c_2$ such that $c_2 - c_1$ is proportional to the element $a$, we will have that
$
E(d^2(V, \{\xi_1 a, \xi_2 a\})) = E(d^2(V, \{c_1^*, c_2^*\})) \leq
E(d^2(V, \{c_1, c_2\}))$, therefore it is possible to determine the principal points of $V$ by considering those of $W_1 = A_1 V$, defining
$c_1^* = \xi_1 a $ and $c_2^* = \xi_2 a $ and then minimizing over $a$.

Therefore, it only remains to obtain $a \in {\cal H}$. Remember
that the operator $A_1$ depends {on} that element
$a$, since it is defined using the normalization of $c_2 - c_1$
which is equal to $a$. To make explicit the dependence, we will
denote it as $A_1^{(a)}$, and also $W_1^{(a)} = A_1^{(a)} V = <a,
V>e_1=Y_1^{(a)}e_1$. Note that $\Sigma_{Y_1^{(a)}} = Var(<a, V>)
= <a, \bGa a>$. Since the principal points will lie in a straight
line with normalized direction $a$, which we are trying to find,
they can be written as $\lambda a$, with $\lambda \in \mathbb{R}$.

Using Lemma \ref{lema:aux:va}, we get that $D_{Y_1^{(a)}}(2) = <a, \bGa a> {D_{\lambda Y_{1}^{(a)}}(2)}/{Var(\lambda Y_{1}^{(a)})}$.
On the other hand, we have that
$$
D_{\lambda Y_{1}^{(a)}}(2) = \min_{\eta_1, \eta_2 \in \mathbb{R}} E\left( \min_{i=1,2} \left\{ |\lambda Y_{1}^{(a)} - \eta_1|^2,
|\lambda Y_{1}^{(a)} - \eta_2|^2 \right\}\right)\le E( |\lambda Y_{1}^{(a)}|^2)=Var\left(\lambda Y_{1}^{(a)}\right)
$$
which implies that ${D_{\lambda Y_{1}^{(a)}}(2)}/{Var(\lambda Y_{1}^{(a)})} < 1$.
Note that by Lemma \ref{lema:aux:va}  we   have that
$$
\frac{D_{\lambda Y_{1}^{(a)}}(2)}{Var(\lambda Y_{1}^{(a)})} = \frac{D_{Y_{1}^{(a)}}(2)}{Var(Y_{1}^{(a)})} $$
and so the ratio does not depend on $\lambda$. Furthermore, we will show that it does not depend on  $a$.

Using that $V$ is elliptic, we get that for any linear and bounded operator $B:{\cal H}\to \mathbb{R}^p$, $\bY=BV$ has an elliptical distribution with parameters $B\mu=0$ and $\bSi=B\bGa B^{*}$. So, its characteristic function can be written as $\varphi_\bY(\by)=\phi(\by\trasp\bSi\by)$ with $\phi$ independent of $B$. In particular, for any $a \in {\cal H}$, we have that $\varphi_{Y_{1}^{(a)}}(b) = \phi(b^2 <a, \bGa a>)=\phi (b^2 Var(Y_{1}^{(a)}))$ which  implies that $Z_a={Y_{1}^{(a)}}/{\sqrt{Var(Y_{1}^{(a)})}}$ has the same distribution for any element $a \in {\cal H}$. Therefore,
\begin{eqnarray*}
\frac{D_{Y_{1}^{(a)}}(2)}{Var(Y_{1}^{(a)})} &=& \frac{\dst\min_{\eta_1, \eta_2 \in \mathbb{R}} E( \min\{ | Y_{1}^{(a)} - \eta_1|^2,
| Y_{1}^{(a)} - \eta_2|^2 \}) } { Var(Y_{1}^{(a)}) }\\
&=& \min_{\eta_1, \eta_2 \in \mathbb{R}} E( \min \left\{ \left( \frac{| Y_{1}^{(a)} - \eta_1|} {\sqrt{Var(Y_{1}^{(a)})}} \right)^2,
\left( \frac{| Y_{1}^{(a)} - \eta_2|} {\sqrt{Var(Y_{1}^{(a)})}} \right)^2 \right\} )\\
&=& \min_{\eta_1^{\star}, \eta_2^{\star} \in \mathbb{R}} E( \min \{ | Z_a - \eta_1^{\star}|^2, | Z_a - \eta_2^{\star}|^2 \} )
\end{eqnarray*}
does not depend on $a$. Hence, we can write ${D_{\lambda Y_{1}^{(a)}}(2)}/{Var(\lambda Y_{1}^{(a)})} = g < 1$ with $g$ independent of $a$ and so $D_{Y_{1}^{(a)}}(2) = g <a, \bGa a>$.

Note that $V =  Y_{1}^{(a)} \, a+ (V -  Y_{1}^{(a)}\, a) = P_{<a>}V + P_{<a>^{\bot}} V$. Then, denoting by $\xi_i^a$ the principal points of $Y_{1}^{(a)}$ and using that that $\|a\| = 1$, we obtain
$$
\|V -  \xi_i^{(a)}\, a\|^2 = \| Y_{1}^{(a)}\, a -   \xi_i^{(a)}\, a + (V -  Y_{1}^{(a)}\, a)\|^2
= \|a \cdot (Y_{1}^{(a)} - \xi_i^a) + P_{<a>^{\bot}} V\|^2 = (Y_{1}^{(a)} - \xi_i^a)^2 + \|P_{<a>^{\bot}} V\|^2.
$$
Taking minimum for $i=1,2$ and then applying expectation, we obtain
\begin{eqnarray*}
E( d^2(V, \{ \xi_1^{(a)}\, a,  \xi_2^{(a)}\, a\}))&=& E ( \min_{i=1,2} \|V - \xi_i^{(a)}\, a\|^2) = D_{Y_{1}^{(a)}}(2) + E (\|P_{<a>^{\bot}} V\|^2)\\
&=&g <a, \bGa a>+ E (\|P_{<a>^{\bot}} V\|^2).
\end{eqnarray*}
Denote $Z = P_{<a>^{\bot}} V$ and let $\phi_j$ be the orthonormal base of $\cal H$ obtained by the eigenfunctions of $\bGa$ related to the eigenvalues $\lambda_1\ge \lambda_2\ge \dots$, then,
\begin{eqnarray*}
E (\|P_{<a>^{\bot}} V\|^2) &=& E(\|Z\|^2) = E(\|V\|^2) - E(\|P_{<a>}V\|^2) = Var(V) - Var(Y_{1}^{(a)})\\
&=& \sum_{n=1}^{\infty}<\bGa \phi_n, \phi_n> - Var(Y_{1}^{(a)}) = \sum_{j\ge 1} \lambda_j- Var(Y_{1}^{(a)})\\
&=&tr(\bGa) - Var(Y_{1}^{(a)}) = tr(\bGa) - <a, \bGa a>.
\end{eqnarray*}
Therefore, we obtain
$$
E( d^2(V, \{ a \cdot \xi_1^a, a \cdot \xi_2^a\})) =  g <a, \bGa a> + tr(\bGa) - <a, \bGa a>
= tr(\bGa)  - (1-g) <a, \bGa a>.
$$
To minimize the left hand side of the above equality it is enough to maximize $<a, \bGa a>$ over the elements $a\in {\cal H}$ with norm equal to 1. Using the compactness of the covariance operator $\bGa$, we obtain the maximum is reached if we choose $a$ as the eigenfunction related to the largest eigenvalue of $\bGa$, concluding the proof. \square

\end{document}